\documentclass[a4paper,11pt]{article}
\usepackage{epsfig}
\newcommand{\sfig}[2]{\parbox{#1 cm}{\psfig{figure=#2}}}
\newcommand{\fsfig}[1]{\parbox{0.25cm}{\psfig{figure=#1, width=0.25cm}}}
\newenvironment{proof}{{\it Proof.\/}}{$\Box$\vskip 0.08in}
\newenvironment{proofl}{{\it Proof of the lemma.\/}}{$\Box$\vskip 0.08in}
\newenvironment{remark}{{\it Remark.\/}}{\vskip 0.08in}
\newenvironment{example}{{\it Example.\/}}{\vskip 0.08in}
\usepackage{amsmath,amssymb,amscd}
\newtheorem{theorem}{Theorem}[section]
\newtheorem{lemma}[theorem]{Lemma}
\newtheorem{corollary}[theorem]{Corollary}
\newcommand{\mc}[1]{\mathcal{#1}}
\newcommand{\mco}[1]{\overline{\mathcal{#1}}}
\begin{document}
\begin{center}
{\huge A spanning tree model for Khovanov homology}
\\[0.2in]
Stephan Wehrli, University of Basel, August 2004 \\[0.2in]
\end{center}
\begin{abstract}
\noindent We use a spanning tree model to prove a result of E.~S.~Lee on the
support of Khovanov homology of alternating knots. \end{abstract}

\ \\
\begin{Large}{\bf Introduction}\end{Large} \bigskip

\noindent In \cite{th}, M.~Thistlethwaite observed that the Kauffman bracket of
a knot diagram is related to the Tutte polynomial of the "black graph" of
the knot diagram. In particular, the Kauffman bracket may be expanded
as a sum over terms corresponding to spanning trees of the "black graph" of the
knot diagram. More than a year ago, the author constructed an analogue of this
expansion for Khovanov homology. The idea to use a spanning tree model for
calculating Khovanov homology was considered by other people
independently.\footnote{
I.~Kofman independently brought up the idea of a spanning tree model for
Khovanov homology at the "Knots in Poland" conference in July 2003,
together with O.~Viro, M.~Polyak, A.~Shumakovitch and L.~Kauffman. I.~Kofman
\cite{ko} gave a talk with the title "Spanning trees and Khovanov homology" at
the "Knots in Washington XVIII" conference in May 2004 (which the author did not
attend). The topic of this paper is very similar to that of I.~Kofman's talk,
yet all results were obtained independently. A paper on the subject by
I.~Kofman, A.~Champanerkar and O.~Viro is in preparation.} In this paper, we use
the spanning tree model for Khovanov homology to give a new proof of a theorem
due to E.~S.~Lee \cite{le1} on the support of the Khovanov homology of
alternating knots. The paper is organized as follows: in Section~\ref{s1}, we
briefly review Thistlethwaite's construction, however without making any
reference to the Tutte polynomial. In Section~\ref{s2}, we describe how the
spanning tree model for the Kauffman bracket leads to a spanning tree model for
Khovanov homology. We give a short proof of a theorem on the behavior of
Khovanov homology under Hopf link addition. This theorem was first proved by
M.~Asaeda and J.~Przytycki in \cite{ap}. In Section~\ref{s3}, we use the results
of Section~\ref{s2} to prove Lee's theorem \cite{le1}. In Section~\ref{s4}, we
discuss a spanning tree model for a homology theory defined in \cite{le2}.

\section{Spanning tree model for the Kauffman bracket}\label{s1}
Unless otherwise stated, link diagrams are assumed to be unoriented and equipped
with a numbering of the crossings. Let $D$ be a link diagram and let
\sfig{0.3}{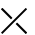} be a crossing of $D$. We may {\it smoothen}
\sfig{0.3}{xing.eps} by replacing it either by \sfig{0.3}{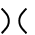} or by
\sfig{0.3}{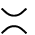}. As in \cite{ba}, we name \sfig{0.3}{0xing.eps} the {\it
$0$-smoothing} and \sfig{0.3}{1xing.eps} the {\it $1$-smoothing} of
\sfig{0.3}{xing.eps}. For a diagram $D'$ obtained from $D$ by smoothening some
of the crossings of $D$, we define $<D|D'>\in\mathbb Z[q,q^{-1}]$ by
$<D|D'>:=(-q)^{r(D,D')}$, where $r(D,D')$ denotes the number of $1$-smoothings
in $D'$. A diagram obtained from $D$ by smoothening all of the crossings of $D$
is called a {\it Kauffman state} of $D$. Let $\mc K(D)$ denote the set of all
Kauffman states of $D$ and $\mc K_k(D)$ the set of all Kauffman states of $D$
consisting of exactly $k$ disjoint circles. Given a link diagram $D$, M.
Khovanov \cite{kh1} assigns a Laurent polynomial $<D>\in\mathbb Z[q,q^{-1}]$ by
the rules
\begin{equation}\label{f1}
<\sfig{0.3}{xing.eps}>=<\sfig{0.3}{0xing.eps}>-q<\sfig{0.3}{1xing.eps}>,
\end{equation}
\begin{equation}\label{f2}
<O^k>=(q+q^{-1})^k.\footnote{\fsfig{xing.eps}, \fsfig{0xing.eps},
\fsfig{1xing.eps} denote any three link diagrams which agree except in a small
disk where they look like \fsfig{xing.eps}, \fsfig{0xing.eps},
\fsfig{1xing.eps}, respectively. $O^k$ denotes a diagram without
crossings which consists of $k$ disjoint circles.}
\end{equation}
$<D>$ is a scaled version of the {\it Kauffman bracket} \cite{ka2}. It is
invariant under Reidemeister moves, up to multiplication with a unit of $\mathbb
Z[q,q^{-1}]$. There is an explicit formula for the Kauffman bracket:
\begin{equation}\label{f3}
<D>=\sum_{D'\in \mc K(D)}<D|D'><D'>.
\end{equation}
To obtain (\ref{f3}) from
(\ref{f1}), we may proceed as follows: First, we apply relation~(\ref{f1}) to
crossing number~$1$ to express $<D>$ as the sum of the two terms on the
right-hand side of (\ref{f1}). Next, we apply relation~(\ref{f1}) to crossing
number~$2$ to express each of these two terms as a sum of two other terms, and
so on. This procedure is visualized in the binary tree in Figure~1:

\ \\
\centerline{\psfig{figure=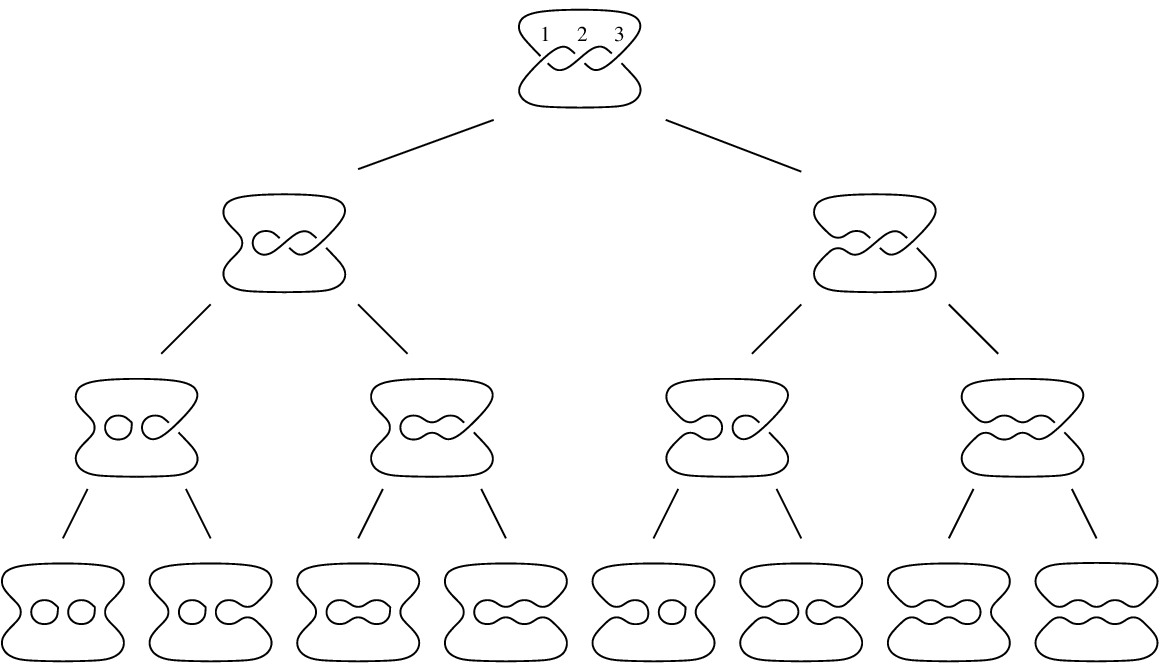,height=4.5cm}}
\begin{center}
Figure 1: Binary tree used to deduce (\ref{f3}) from (\ref{f1}).
\end{center}
The diagrams sitting at the leaves of the tree are the
Kauffman states of $D$, whence (\ref{f3}) follows.

Unfortunately, the complexity of (\ref{f3})
grows exponentially in the complexity of $D$. In case $D$ is connected, we get a
more efficient formula by modifying the above procedure as follows: each time
before applying relation~(\ref{f1}), we check the connectivity of the diagrams
on the right-hand side. We rewrite the term on the left-hand side as the sum of
the terms on the right-hand side only if both diagrams on the right-hand side
are connected. The modified procedure is visualized in the binary tree in
Figure~2:

\ \\
\centerline{\psfig{figure=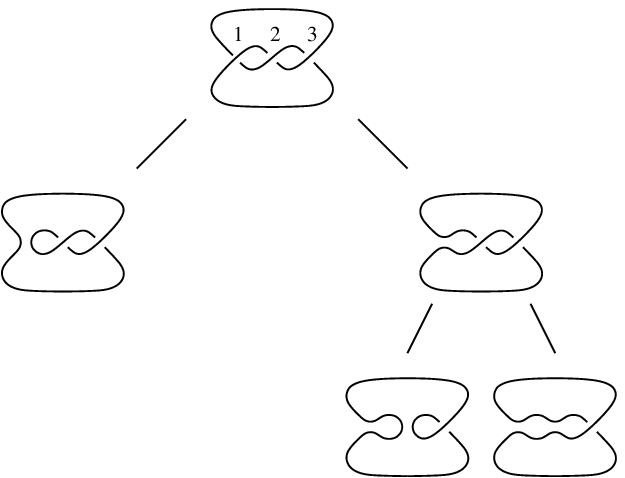,height=3cm}}
\begin{center}
Figure 2: Binary tree used to deduce (\ref{f4}) from (\ref{f1}).
\end{center}
We immediately obtain:

\begin{equation}\label{f4}
<D>=\sum_{D'\in\mc T(D)}<D|D'><D'>,
\end{equation}
where $\mc T(D)$ denotes the set of all diagrams sitting at the leaves of the
tree in Figure~2. (Note that $\mc T(D)$ may depend on the numbering of the
crossings of $D$.) Let $D'$ be an element of $\mc T(D)$. By construction $D'$ is
connected and every crossing of $D'$ is {\it splitting} (i.e. connects two
otherwise disconnected parts of $D'$). A diagram with these properties will be
called {\it R1-trivial} because it can be trivialized using Reidemeister move~1
only.

It is easy to see that the Kauffman bracket behaves as follows under
Reidemeister move~1:
\begin{equation}\label{f5}
<\sfig{0.35}{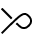}>=q^{-1}<\sfig{0.15}{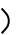}>,\qquad
<\sfig{0.35}{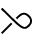}>=-q^2<\sfig{0.15}{notwist.eps}>.
\end{equation}
From (\ref{f5}) we get an explicit formula for the Kauffman
bracket of an R1-trivial diagram: Let $D'$ be an R1-trivial diagram and assume
that an orientation of $D'$ is given. Let $x(D')$ denote the number of negative
(\sfig{0.3}{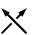}) crossings of $D'$ and $y(D')$ the number of positive
(\sfig{0.3}{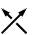}) crossings of $D'$. Then
\begin{equation}\label{f6}
<D'>=(-1)^{x(D')}q^{2x(D')-y(D')}(q+q^{-1}).
\end{equation}

Since $\mc K_1(D)=\bigsqcup_{D'\in \mc T(D)}\mc K_1(D')$ and since $\#\mc
K_1(D')=1$ for $D'$ R1-trivial, the elements of $\mc K_1(D)$
correspond bijectively to the elements of $\mc T(D)$. Indeed, for $S\in\mc
K_1(D)$ let $D_S$ be the unique element of $\mc T(D)$ having $S$ among its
Kauffman states. We may rewrite (\ref{f4}) as
\begin{equation}\label{f7}
<D>=\sum_{S\in \mc K_1(D)}<D|D_S><D_S>.
\end{equation}

Assume that the regions of the knot projection underlying the knot diagram $D$
are colored black and white in a chessboard pattern, i.e. such that any two
regions which share an edge have opposite colors and such that the unbounded
region is colored white. (There is exactly one such coloring. It is obtained by
coloring the bounded and the unbounded region of any $S\in\mc K_1(D)$ black and
white, respectively.) A smoothing of a crossing of $D$ will be called a {\it
black smoothing} or a {\it white smoothing} depending on whether it connects two
black regions or two white regions of $D$. The {\it black graph} of $D$ is the
graph whose vertices correspond to the black regions of $D$ and which has an
edge connecting two vertices for each crossing where the corresponding black
regions touch. There is a bijection between spanning trees of the black graph of
$D$ and elements of $\mc K_1(D)$ defined as follows: to a spanning tree assign
an element of $\mc K_1(D)$ by choosing the black smoothing for precisely those
crossings of $D$ which correspond to edges of the spanning tree. We call
(\ref{f7}) a {\it spanning tree model} for the Kauffman
bracket.\footnote{Although our approach is different, formula~(\ref{f7}) is
essentially the same as the formula used to define $\Gamma_G$ in \cite{th}. Note
that both formulae depend on a numbering of the crossings. To go from (\ref{f7})
to the fomula in \cite{th}, one has to reverse the order of the crossings.}

Since the number of edges in a tree is one less than the number of vertices, the
number of black smoothings in an element of $\mc K_1(D)$ is one less than the
number of black regions of $D$. In particular, the number of black smoothings is
the same for all $S\in\mc K_1(D)$. This fact will be important, so let us give
another proof of it: by the Clock Theorem \cite{ka1}, any two elements
of $\mc K_1(D)$ are related by a finite sequence of state transpositions (see
Figure~3) and state transpositions do not change the number of black smoothings.

\ \\
\centerline{
\begin{tabular}{c@{\qquad}c@{\qquad}c}
\psfig{figure=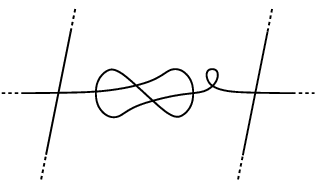,width=3.7cm}&
\psfig{figure=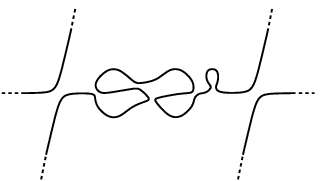,width=3.7cm}&
\psfig{figure=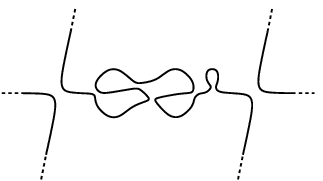,width=3.7cm}
\end{tabular}
}
\begin{center}
Figure 3: A knot projection and two smoothings related by
a state transposition.
\end{center}
The fact that the number of black smoothings in $S\in\mc K_1(D)$ is independent
of $S$ may also be shown by induction on the number of crossings of $D$.


\section{Spanning tree model for Khovanov homology}\label{s2}

In \cite{kh1}, M. Khovanov assigns to a link diagram $D$ a bigraded complex
$\mco C(D)$ with differential $d$ of bidegree $(1,0)$. Let $\mco H(D)$ denote
the homology of $\mco C(D)$. We call $\mco C(D)$ the {\it Khovanov complex} of
$D$ and $\mco H(D)$ the {\it Khovanov homology}. The Khovanov complex is a link
invariant when considered up to shifts of the gradings and up to a chain
equivalence which preserves the secondary grading
(and also the primary grading, but this is already contained in the
definition of a chain equivalence). M. Khovanov's construction may be viewed
as a categorification of the Kauffman bracket. Indeed, the Kauffman
bracket is the "graded Euler characteristic" of the Khovanov complex (see
\cite{kh1}). In this section, we construct an analogue of (\ref{f7}) for the
Khovanov complex. We use the following notations and conventions: all modules
and complexes are assumed to be bigraded. Isomorphisms between bigraded objects
are assumed to preserve the gradings. Direct sums are assumed to be compatible
with the gradings. If $M$ is a bigraded object, we denote by $M^{i,j}$ its
homogeneous component of bidegree $(i,j)$. For $m,n\in\mathbb Z$, $[m]$ and
$\{n\}$ denote the operators on bigraded objects defined by
$(M[m]\{n\})^{i,j}:=M^{i-m,j-n}$. Unless otherwise stated, we work with
coefficients in $\mathbb Z$.

Instead of giving the precise definition of $\mco C(D)$, we list the
properties of $\mco C(D)$ which are relevant to our discussion. On the level of
modules,
\begin{equation}\label{f8}
\mco C(D)=\bigoplus_{D'\in\mc K(D)} \mco C(D')[r(D,D')]\{r(D,D')\}.
\end{equation}
In particular, $\mco C(\sfig{0.3}{xing.eps})$ as a module is canonically
isomorphic to $\mco C(\sfig{0.3}{0xing.eps})\oplus\mco
C(\sfig{0.3}{1xing.eps})[1]\{1\}$.
On the level of complexes, the modification \sfig{0.3}{0xing.eps}
$\rightsquigarrow$ \sfig{0.3}{1xing.eps} induces a chain transformation
\begin{equation}\label{f9}
w:\mco C(\sfig{0.3}{0xing.eps})\longrightarrow \mco
C(\sfig{0.3}{1xing.eps})\{1\}
\end{equation}
which preserves the secondary grading and such that $\mco
C(\sfig{0.3}{xing.eps})$ is canonically isomorphic to the mapping
cone\footnote{Let $w$ be a chain transformation from a complex $C_0$ with
differential $d_0$ to a complex $C_1$ with differential $d_1$. The {\it mapping
cone} of $w$ is the complex $C$ with differential $d$ defined as follows: As a
module, $C=C_0\oplus C_1[1]$. The restriction of $d$
to $C_0$ is $d_0+w$. The restriction of $d$ to $C_1[1]$ is
$-d_1$.} of $w$. The Khovanov complex of a diagram without crossings has trivial
differential. As a module, it is given by
\begin{equation}\label{f10} \mco
C(O^k)=\mc A^{\otimes k} \end{equation}
where  $\mc A$ is the bigraded module defined by $\mc A^{0,-1}=\mc
A^{0,1}=\mathbb Z$ and $\mc A^{i,j}=0$ for $(i,j)\neq(0,\pm 1)$. The gradings
are additive under tensor multiplication. By \cite[Section 5]{kh1}, the behavior
of the Khovanov complex under Reidemeister move~1 is as follows:
\begin{equation}\label{f11} \mco C(\sfig{0.35}{ltwist.eps}) \cong \mco
C(\sfig{0.15}{notwist.eps})\{-1\} \oplus B_1,\qquad\mco
C(\sfig{0.35}{rtwist.eps}) \cong \mco C(\sfig{0.15}{notwist.eps})[1]\{2\} \oplus
B_2 \end{equation}
for contractible\footnote{A complex $B$ with
differential $d$ of bidegree $(1,0)$ is called {\it contractible} if it is chain
equivalent to the trivial complex, i.e. if there exists a graded module
endomorphism $p$ of $B$ of bidegree $(-1,0)$ such that $p\circ d+d\circ p=id_B$.
Note that a contractible complex has trivial homology.} complexes $B_1$, $B_2$.

We are now ready to discuss how (\ref{f7}) transfers to the Khovanov complex.
Consider the Khovanov complex of a diagram sitting at an internal node in the
binary tree of Figure~2. It is canonically isomorphic to the mapping cone of a
chain transformation between the Khovanov complexes of the two diagrams sitting
right below the given diagram in the binary tree. As a module, it is equal to
the direct sum of the modules underlying these two complexes. (Actually, there
are also some shifts of the gradings.) Hence on the level of modules
\begin{equation}\label{f12} \mco C(D)=\bigoplus_{S\in \mc K_1(D)}\mco
C(D_S)[r(D,D_S)]\{r(D,D_S)\}. \end{equation}
By (\ref{f11}), the complexes $\mco C(D_S)$ admit decompositions
\begin{equation}\label{f13}
\mco C(D_S)\cong\mc A[x(D_S)]\{2x(D_S)-y(D_S)\}\oplus B(D_S),
\end{equation} for contractible complexes $B(D_S)$. To obtain a decomposition of
the complex $\mco C(D)$, we need the following lemma which asserts that forming
the direct sum with a contractible complex "commutes", up to isomorphism, with
the mapping cone construction:

\begin{lemma}\label{l2.1} Let $C_0$ and $C_1$
be complexes with $C_i=A_i\oplus B_i$ for complexes $A_i$,
$B_i$ with $B_i$ contractible. Let $w:C_0\rightarrow C_1$ be a
grading-preserving chain transformation and let $w_{AA}:A_0\rightarrow A_1$
denote $w$ composed with the obvious projection and inclusion. Let $A$, $B$ and
$C$ be the mapping cone of $w_{AA}$, the (contractible) complex $B_0\oplus
B_1[1]$ and the mapping cone of $w$, respectively. Then $C\cong A\oplus B$.
\end{lemma} \begin{proof} Define
$w_{AB}:B_0\rightarrow A_1$, $w_{BA}:A_0\rightarrow B_1$ and
$w_{BB}:B_0\rightarrow B_1$ in the same way as $w_{AA}$. Since $B_i$ is
contractible, there exists $p_i$ with $p_i\circ d_i + d_i\circ p_i = id_{B_i}$.
Writing elements of both $C$ and $A\oplus B$ in the form $(a_0,b_0,a_1,b_1)$,
$a_i \in A_i$, $b_i \in B_i$, we define $f:C\rightarrow A\oplus B$ by
$f(a_0,b_0,a_1,b_1):=(a_0,b_0,a_1-(w_{AB}\circ p_0)b,b_1-(w_{BB}\circ
p_0)b-(p_1\circ w_{BA})a)$. Direct calculation shows that $f$ is an isomorphism
of complexes.
\end{proof}
In view of (\ref{f12}), (\ref{f13}) and Lemma~\ref{l2.1}, it is easy to see that
the complex $\mco C(D)$ is isomorphic to a direct sum of two complexes, $A$ and
$B$ say, where $B$ is the direct sum of the $B(D_S)$ (up to shifts of the
gradings) and $A$ as a module is given by
\begin{equation}\label{f14}
A=\bigoplus_{S\in \mc K_1(D)}\mc
A[x(D_S)]\{2x(D_S)-y(D_S)\}[r(D,D_S)]\{r(D,D_S)\}. \end{equation}
Writing $w(D_S)$ for $x(D_S)-y(D_S)$ and observing that
$r(D,D_S)=r(D,S)-r(D_S,S)=r(D,S)-y(D_S)$, we obtain the following theorem:

\begin{theorem}\label{t2.2}
Let $D$ be a connected link diagram. Then there is a decomposition $\mco
C(D)\cong A\oplus B$, where $B$ is contractible and $A$ as a module is given by
\begin{equation}\label{f15}
A=\bigoplus_{S\in \mc K_1(D)}\mc
A[w(D_S)]\{2w(D_S)\}[r(D,S)]\{r(D,S)\}.
\end{equation}
We call the decomposition $\mco C(D)\cong A\oplus B$ together with the above
formula a {\normalfont spanning tree model} for Khovanov homology.
\end{theorem}
Although the above constructions are completely explicit, it may be difficult to
compute the differential of $A$ in practice. However,
we have the estimates
$\mbox{dim}_{\mathbb Q}(\mco H(D)\otimes\mathbb Q)\leq\mbox{rank}(A)$ and
$\mbox{dim}_{\mathbb Q}(\mco H(D)^{i,j}\otimes\mathbb
Q)\leq\mbox{rank}(A^{i,j})$, which together with (\ref{f15}) imply the following
theorem:
\begin{theorem}\label{t2.3}
Let $D$ be a connected link diagram. Then
\begin{equation}
\mbox{\normalfont dim}_{\mathbb Q}(\mco
H(D)\otimes\mathbb Q)\leq 2(\#\mc K_1(D)).
\end{equation}
Moreover, the dimension of $\mco
H(D)^{i,j}\otimes\mathbb Q$ is bounded from above by the number of
$S\in\mc K_1(D)$ with $w(D_S)+r(D,S)=i$ and $r(D,S)=2i-j\pm 1$.
\end{theorem}
Theorem~\ref{t2.3} shows
that the dimension of $\mco H(D)\otimes\mathbb Q$ tends to be much smaller than
the rank of $\mco C(D)$, a fact observed experimentally by D.~Bar-Natan
\cite{ba}.

\medskip
\noindent\begin{remark} In \cite{kh2}, $H^1:=\mc A\{-1\}$ is endowed
with the structure of a $\mathbb Z_-$-graded\footnote{Our grading on $H^1$ is
opposite to the grading used in \cite{kh2}.} commutative ring. In
particular, $\mc A$ is an $H^1$-module and $\mc A^{0,-1}\subset\mc A$ an
$H^1$-submodule. When a distinguished point on the knot projection (which is not
a double point) is chosen, $\mco C(D)$ becomes a complex of $H^1$-modules. The
isomorphisms in (\ref{f11}) are isomorphisms of $H^1$-modules, provided
Reidemeister move~$1$ is performed away from the distinguished point. Hence we
may assume that the decomposition of $\mco C(D)$ given in Theorem~\ref{t2.2} is
compatible with the $H^1$-module structure. By (\ref{f15}), $A$ is generated as
an $H^1$-module by elements corresponding to elements of $\mc K_1(D)$, i.e. to
spanning trees of the black graph of $D$. Similarly, the {\it reduced Khovanov
complexes} $\mco C(D)\otimes_{H^1}\mc A^{0,-1}$ and $\mco C(D)\otimes_{H^1}(\mc
A/\mc A^{0,-1})$ are chain equivalent to complexes which are generated as
$\mathbb Z$-modules by elements corresponding to spanning trees of the black
graph of $D$. It is interesting to compare this with knot Floer homology (see
\cite{os}). \end{remark}

\noindent \begin{remark} By \cite{kh2}, there is an isomorphism $\mco C(D_1\#
D_2)\cong\mco C(D_1)\otimes_{H^1}\mco C(D_2)\{-1\}$. (Here, the $H^1$-module
structures on $\mco C(D_1)$ and $\mco C(D_2)$ are defined as in the previous
remark, by choosing distinguished points on $D_1$ and $D_2$ close to the
connected sum point.) Let $D_2$ be a standard diagram of the Hopf link. By
Theorem \ref{t2.2}, $\mco C(D_2)\cong A\oplus B$ where $B$ is contractible and
$A$ as an $H^1$-module is equal to $\mc A[-1]\{-2\}\oplus\mc A[1]\{2\}$. The
differential of $A$ has bidegree $(1,0)$, whence it must be zero. We conclude
that $\mco C(D_1\#D_2)$ is chain equivalent to $\mco C(D_1)[-1]\{-2\}\oplus\mco
C(D_1)[1]\{2\}$. Consequently, $\mco H(D_1\#D_2)$ is isomorphic to
$\mco H(D_1)[-1]\{-2\}\oplus\mco H(D_1)[1]\{2\}$. The latter was first proved
by M.~Asaeda and J.~Przytycki \cite{ap}, who thus confirmed a conjecture of
A.~Shumakovitch \cite{sh1}. A special case was already considered in
\cite[Section 4]{we}. \end{remark}


\section{Support of Khovanov homology for alternating knots}\label{s3}

The theorems in this section were conjectured by D.~Bar-Natan, S.~Garoufalidis
and M.~Khovanov \cite{ba} and proved by E.~S.~Lee \cite{le1}. We give new proofs
using the spanning tree model. For short proofs, see also \cite{ap}.

A knot diagram is {\it alternating} if one alternately over- and undercrosses
other strands as one goes along the knot in that diagram. A knot is called {\it
alternating} if it possesses an alternating diagram.

Let $D$ be an alternating knot diagram. Then either the black
smoothing coincides with the $1$-smoothing for all crossings or the white
smoothing coincides with the $1$-smoothing for all crossings. By
Section~\ref{s1}, the numbers of black and white smoothings in $S\in\mc K_1(D)$
are independent of $S$. Therefore $r(D,S)$ is independent of $S\in\mc K_1(D)$.
Let $n_1(D):=r(D,S)$ for any $S\in\mc K_1(D)$. From Theorem~\ref{t2.2} we get:
\begin{theorem}\label{t3.1}
Let $D$ be an alternating knot diagram. Then $\mco
H(D)^{i,j}=0$ unless $(i,j)$ lies on one of the two lines $j=2i-n_1(D)\pm 1$.
Moreover, $\mco H(D)^{i,j}$ is torsion free unless $j=2i-n_1(D)-1$.
\end{theorem}
The statement about the torsion follows from the fact that the differential
of $A$ has bidegree $(1,0)$. We also obtain:
\begin{theorem}\label{t3.2}
Let $D$ be an alternating knot diagram. Let $i_-$ and $i_+$ be the smallest and
largest primary degree in which $\mco H(D)$ is non-zero. Then $\mco
H(D)^{i_-,j_-}$ and $\mco H(D)^{i_+,j_+}$ are non-zero for $(i_-,j_-)$ on the
lower line and $(i_+,j_+)$ on the upper line (i.e. $j_-=2i_--n_1(D)-1$ and
$j_+=2i_+-n_1(D)+1$). \end{theorem}
The following theorem corresponds to parts (i) and (iii) of Theorem~1 of
\cite{th}. Our proof will be related to the proofs given in
\cite{th}. For a different proof of a similar statement, see \cite[Section
7.7]{kh1}.

\begin{theorem}\label{t3.3}
Let $D$ be an alternating knot diagram. Assume that no crossing of $D$ is
splitting. Then $\mco H(D)^{i_-,j_-}=\mco H(D)^{i_+,j_+}=\mathbb Z$ and $i_-=0$
and $i_+$ is equal to the number of crossings of $D$.  \end{theorem}
Theorem~\ref{t3.3} is a consequence of Theorem~\ref{t2.2} and the following
lemma, in which $n_0(D)$ denotes the number of $0$-smoothings in any $S\in\mc
K_1(D)$:

\begin{lemma}
Let $D$ be an alternating knot diagram. Assume that no crossing of $D$ is
splitting. For any $S\in\mc K_1(D)$ there is a numbering of the crossings of
$D$ such that $w(D_S)=-n_1(D)$ and $w(D_{S'})>-n_1(D)$ for $S'\neq S$.
Likewise, there is a numbering of the crossings of $D$ such that $w(D_S)=n_0(D)$
and $w(D_{S'})<n_0(D)$ for $S'\neq S$.
\end{lemma}
\begin{proofl}
Let $S\in\mc K_1(D)$. Let us number the crossings of $D$ in such
a way that the crossings which are $0$-smoothings in $S$ precede those which are
$1$-smoothings in $S$. Since either the black smoothing coincides
with the $1$-smoothing for all crossings or the white smoothing
coincides with the $1$-smoothing for all crossings and since no crossing
of $D$ is splitting, the construction of $D_S$ implies that
the crossings of $D$ which are not smoothened in $D_S$ are precisely those which
are $1$-smoothings in $S$. Moreover, a look at local orientations shows that
these crossings have to be positive with respect to any orientation of $D_S$,
for otherwise $S$ could not be connected. So $w(D_S)=-y(D_S)=-n_1(D)$. Now let
us consider $S'\in\mc K_1(D)$ with $S'\neq S$. The first crossing of $D$ where
$S$ and $S'$ differ has to be a $0$-smoothing in $S$ and a $1$-smoothing in
$S'$. Moreover, it must be smoothened in $D_{S'}$ because it is smoothened in
$D_S$. We conclude $y(D_{S'})<y(D_S)$, whence $w(D_{S'})>-y(D_S)=-n_1(D)$.
Analogously, if we number the crossings of $D$ in such a way that the
crossings which are $1$-smoothings in $S$ precede those which are $0$-smoothings
in $S$, we have $w(D_S)=n_0(D)$ and $w(D_{S'})<n_0(D)$ for $S'\neq S$.
\end{proofl}
Note that the difference $i_+(D)-i_-(D)$ is a knot invariant. By
(\ref{f8}), it cannot exceed the number of crossings of
$D$. From Theorem~\ref{t3.3} we get the following corollary which coincides with
Corollary~1 of \cite{th}:

\begin{corollary}
If a knot possesses an alternating diagram with $m$ crossings, all of
which are non-splitting, then the knot does not admit a diagram with fewer than
$m$ crossings. \end{corollary}

\section{Spanning tree model and Lee's differential}\label{s4}
In \cite{le2}, E.~S.~Lee defined a differential $\Phi$ of bidegree $(1,4)$ on
$\mco C(D)$ which anticommutes with the differential
$d$. It follows that $(d+\Phi)^2=0$, whence $d':=d+\Phi$ may be considered as a
differential on $\mco C(D)$. Note that $d'$ does not decrease the secondary
degree. Otherwise stated, $d'$ respects the filtration
defined as follows: an element of $\mco C(D)$ has filtration at least $j$ if and
only if it is a sum of homogeneous elements of secondary degree at least $j$.
Let $\mco C'(D)$ denote the complex $\mco C(D)$ with differential $d'$ instead
of $d$.
\begin{theorem}\label{t4.1}
The complex $\mco C'(D)$ is invariant under Reidemeister moves, up to shift of
the primary grading and up to filtered chain equivalence. In particular, there
are decompositions
\begin{equation}
\mco C'(\sfig{0.35}{ltwist.eps}) \cong \mco
C'(\sfig{0.15}{notwist.eps})\{-1\} \oplus B'_1,\qquad \mco
C'(\sfig{0.35}{rtwist.eps}) \cong \mco C'(\sfig{0.15}{notwist.eps})[1]\{2\}
\oplus B'_2
\end{equation}
for contractible filtered complexes $B'_1$ and $B'_2$.
These decompositions respect the filtration (meaning that
the associated inclusion and projection maps respect the filtration).
\end{theorem}
To prove Theorem~\ref{t4.1}, we may adopt the proofs of the corresponding
statements for the Khovanov complex given in \cite[Section 5]{kh1}. Let us
explain this in more detail. Let $\mc M$ denote the category which has
closed $1$-manifolds as objects and cobordisms as morphisms (and disjoint union
as tensor product). The definition of the Khovanov complex given in \cite{kh1}
involves a monoidal functor $\mc F$ from $\mc M$ to the category which has
graded $\mathbb Z$-modules as objects and graded $\mathbb Z$-module
homomorphisms as morphisms. On objects, $\mc F$ is given by $\mc F(O^k)=\mc
A^{\otimes k}$. $\mc F$ satisfies $\mc F(S^0_0)=0$ and
\begin{equation}\label{f18} \mc A\otimes\mc A=\mc F(S^1_1\sqcup S^1_0)\mc
A\oplus\mc F(S^2_1)\mc A, \end{equation}
where $S^l_k$ denotes the cobordism
from $O^k$ to $O^l$ which is a $2$-sphere with $k+l$ disks removed. Note that
the inclusion and projection maps associated with decomposition~(\ref{f18}) may
be described in terms of the functor $\mc F$. (For example, the projection onto
$\mc F(S^1_1\sqcup S^1_0)\mc A$ is given by $\mc F(S^1_2)\circ(\mc F(S^1_1\sqcup
S^1_1)-\mc F(S^2_1)\circ\mc F(S^1_1\sqcup S^0_1))$. To verify this,
use $\mc F(S^0_0)=0$ and the functoriality of $\mc F$.) In \cite{ra}, J.~Rasmussen
remarked that $\mco C'(D)$ may be defined by using a monoidal functor
$\mc F'$ from $\mc M$ to the category which has filtered $\mathbb Z$-modules as
objects and filtered $\mathbb Z$-module homomorphisms as morphisms. $\mc F'$
satisfies $\mc F'(O^k)=\mc A^{\otimes k}$ and $\mc F'(S^0_0)=0$. It is
straightforward to see that (\ref{f18}) remains satisfied when $\mc F$ is
replaced by $\mc F'$. Now Theorem~\ref{t4.1} may be established by replacing
$\mc F$ by $\mc F'$ in the proofs of \cite[Section 5]{kh1}. Note that these
proofs rely only on (\ref{f18}) and on the functoriality of $\mc F$.

From Theorem~\ref{t4.1} and the definition of
$\mco C'(D)$ it is clear that the properties of the Khovanov complex used to
deduce Theorem~\ref{t2.2} have analogues for $\mco C'(D)$. Therefore, we obtain:
\begin{theorem}
Let $D$ be a connected link diagram. Then $\mco C'(D)\cong A\oplus B$ where
$B$ is a contractible filtered complex and $A$ as a filtered module is given by
(\ref{f15}). The decomposition $\mco C'(D)\cong A\oplus B$ respects the
filtration.
\end{theorem}

The remarks at the end of Section~\ref{s2} remain true for $\mco C'(D)$ if the
multiplication on $H^1=\mc A\{-1\}$ is replaced by the multiplication induced by
$\mc F'(S^1_2):\mc A\otimes\mc A\rightarrow\mc A$ and $\mc A^{0,-1}\subset\mc
A$ is replaced by $\mathbb Za\subset\mc A$, for $a\in\mc A$ defined as in
\cite{le2}.
Let $\mco H'(D)$ denote the bigraded module given by $\mco
H'(D)^{i,j}:=F^{i,j}/F^{i,j-1}$, where $F^{i,j}$ is the submodule of the
homology of $\mco C'(D)$ consisting of all homology classes which have
representatives of primary degree $i$ and filtration at least $j$. The first
part of Theorem~\ref{t3.1} remains true if $\mco H(D)$ is replaced by $\mco
H'(D)$. In general, the statement about the torsion is not true for $\mco
H'(D)$.
\medskip \ \\
\begin{example} Let $D$ be a standard diagram of the left
handed trefoil. Then $\mco C'(D)\cong A\oplus B$ for a contractible complex $B$
and $A=\mc A[0]\{-1\}\oplus\mc A[2]\{3\}\oplus\mc A[3]\{5\}$, as a module. The
differential of $A$ is zero on $\mc A[0]\{-1\}$ and maps $1,X\in\mc A[2]\{3\}$
to $2X,2\cdot 1\in\mc A[3]\{5\}$, respectively. Whence $\mco H'(D)\cong\mc
A[0]\{-1\}\oplus (\mc A/2\mc A)[3]\{5\}$. In particular, $\mco
H'(D)^{3,6}\cong\mathbb Z/2\mathbb Z$, despite the fact that $(3,6)$ lies
on the upper of the two lines mentioned in Theorem~\ref{t3.1}.
\end{example}
\noindent The example also shows that Theorem~\ref{t3.3} is not true for $\mco
H'(D)$. In this context, let us recall some facts from
\cite{le2} and \cite{ra}:  If $D$ is a knot diagram, all elements of $\mco
H'(D)\otimes\mathbb Q$ have the same primary degree. For any diagram $D$, the
filtration on $\mco C'(D)$ induces a spectral sequence whose $E_2$-term is $\mco
H(D)$ and which converges to $\mco H'(D)$. Therefore,
$\operatorname{dim}_{\mathbb Q}(\mco H'(D)\otimes\mathbb
Q)\leq\operatorname{dim}_{\mathbb Q}(\mco H(D)\otimes\mathbb Q)$. If $D$ is a
diagram of a $k$-component link, $\operatorname{dim}_{\mathbb Q}(\mco
H'(D)\otimes\mathbb Q)=2^k$. \medskip \ \\
\begin{remark} There is a decomposition of $\mco
C'(D)\otimes\mathbb Q$ (which does not respect the filtration) into a
contractible complex and a complex of dimension $2^k$ with trivial differential.
Let us briefly describe this decomposition. Let $\{a,b\}\subset \mc A\otimes
\mathbb Q$ be the basis of $\mc A\otimes \mathbb Q$ defined in \cite{le2}. We
call $a$ and $b$ {\it colors}. A {\it colored Kauffman state} is a Kauffman
state $D'$ together with an assignment of $a$ or $b$ to every circle of $D'$. By
(\ref{f8}) and (\ref{f10}), $\mco C'(D)\otimes\mathbb Q$ is spanned by vectors
corresponding to colored Kauffman states of $D$.\footnote{This description of
$\mco C'(D)\otimes\mathbb Q$ is similar to the description of $\mco C(D)$ given
in \cite{vi}.} A {\it coloring} of a diagram $D$ is an assignment of $a$ or $b$
to every edge of $D$. A coloring of $D$ is called {\it admissible} if for every
crossing of $D$ either all four edges touching at the crossing have the same
color or two neighboring edges (for example the lower left edge and the lower
right edge in \sfig{0.3}{xing.eps}) have color $a$ and the other two edges have
color $b$. If $c$ is an admissible coloring of $D$, we may smoothen every
two-color crossing of $D$ in the way consistent with the coloring to obtain a
colored diagram $D_c$ in which every component is colored consistently. Let
$V(D_c)$ denote the subspace of $\mco C'(D)\otimes\mathbb Q$ spanned by
all colored Kauffman states of $D$ whose colorings agree with the coloring of
$D_c$. The structure of $d'$ implies that $V(D_c)$ is a subcomplex of
$\mco C'(D)\otimes\mathbb Q$, whence we get a decomposition of $\mco
C'(D)\otimes\mathbb Q$ (which does not respect the filtration):
\begin{equation}
\mco C'(D)\otimes\mathbb Q=\bigoplus_{c\,\,\mbox{\scriptsize admissible}}
V(D_c). \end{equation}
The subcomplexes $V(D_c)$ are easy to understand: if $D_c$ has at least
one crossing, $V(D_c)$ is isomorphic to the mapping cone of an isomorphism
and hence conctractible. If $D_c$ has no crossings, $V(D_c)$ is
one-dimensional and hence has trivial differential. The admissible colorings $c$
of $D$ for which $D_c$ has no crossings correspond bijectively to the possible
orientations of $D$ (see \cite{le2} and \cite{ra}). Hence there are $2^k$ such
colorings. \end{remark}
\ \\
\begin{Large}{\bf Acknowledgements}\end{Large} \bigskip

\noindent The author was partially sponsored by the Swiss National Science
Foundation. The author wishes to thank Norbert A'Campo and
Sebastian Baader for many stimulating conversations.

\noindent
{\it E-mail address:} \texttt{wehrli@math-lab.unibas.ch}

\begin{thebibliography}{000000}
\bibitem[AP]{ap} Marta M. Asaeda, J\'ozef H. Przytycki, {\it Khovanov
homology: torsion and thickness}, arXiv:math.GT/0402402, 2004.
\bibitem[BN]{ba} Dror Bar-Natan, {\it On Khovanov's Categorification of the
Jones  polynomial}, Algebraic and Geometric Topology 2, 2002, 337-370,
arXiv:math.QA/0201043 .
\bibitem[Ka-1]{ka1} Louis H. Kauffman, {\it Formal Knot
Theory}, Mathematical Notes 30, Princeton University Press, 1983.
\bibitem[Ka-2]{ka2} Louis H. Kauffman, {\it State models and the Jones
polynomial}, Topology 26, 1987, 395-407.
\bibitem[Kh-1]{kh1} Mikhail
Khovanov, {\it A categorification of the Jones polynomial}, Duke Math. J.
101, 2000, no. 3, 359-426, arXiv:math.QA/9908171, 2001.
\bibitem[Kh-2]{kh2} Mikhail Khovanov, {\it A functor-valued invariant of
tangles}, arXiv:math.QA/0103190, 2001.
\bibitem[Ko]{ko} Ilya S. Kofman, Private communication; an abstract of the
talk on spanning trees and Khovanov homology at the "Knots in Washington XVIII"
conference in May 2004 can be found at http://at.yorku.ca/c/a/n/v/19.htm .
\bibitem[ML]{ma} Saunders Mac Lane, {\it Homology}, Die Grundlehren der mathema-
tischen Wissenschaften Band 114, Springer, 1963.
\bibitem[Lee-1]{le1} Eun Soo Lee, {\it The support of Khovanov's
invariants for alternating knots}, arxiv:math.GT/0201105, 2002.
\bibitem[Lee-2]{le2} Eun Soo Lee, {\it On Khovanov invariant for alternating
links}, arXiv:math.GT/0210213, 2002.
\bibitem[OSz]{os} Peter Ozsv\'ath, Zolt\'an
Szab\'o, {\it Knot Floer homology, genus bounds, and mutation},
arXiv:math.GT/0303225, 2003.
\bibitem[Ra]{ra} Jacob Rasmussen, {\it Khovanov
homology and the slice genus}, arXiv.math.GT/0402131, 2004.
\bibitem[Sh-1]{sh1} Alexander Shumakovitch, Private communication.
\bibitem[Sh-2]{sh2} Alexander Shumakovitch, {\it Torsion of the Khovanov
homology}, arXiv:math.GT/0405474, 2004.
\bibitem[Th]{th} Morwen B. Thistlethwaite , {\it A spanning tree expansion of
the Jones polynomial}, Topology 26, 1987, 297-309.
\bibitem[Vi]{vi} Oleg Viro, {\it Remarks on the definition of Khovanov
Homology}, arXiv:math.GT/0202199, 2002.
\bibitem[We]{we} Stephan M. Wehrli, {\it Khovanov Homology and Conway Mutation},
arXiv:math.GT/0301312, 2003.
\end{thebibliography}
\end{document}